\documentclass[10pt,english]{article}
\usepackage[T1]{fontenc}
\usepackage[latin9]{inputenc}
\usepackage[paperwidth=6.5in,paperheight=9.2in]{geometry}
\usepackage{amstext}
\usepackage{amssymb}
\usepackage{amsmath}
\usepackage{graphicx}
\usepackage{setspace}



\usepackage{babel}
\begin{document}

\title{Is the spiral effect psychological?}
\author{Bernhard Klaassen}
\date{(To appear in \emph{Elemente der Mathematik} 2022)}
\maketitle

\section{Introduction and definitions}
In 2017 \cite{klaassen} a definition of spiral tilings was given, thereby
answering a question posed by Gr\"unbaum and Shephard in the late 1970s.
The author had the pleasure to discuss the topic via e-mail with Branko
Gr\"unbaum in his 87th year. During this correspondence the question
arose whether a spiral structure (given a certain definition of it)
could be recognized automatically or whether ``to some
extent, at least, the spiral effect is psychological\textquotedblright ,
as Gr\"unbaum and Shephard had conjectured in 1987 (see exercise section of chapter 9.5 in \cite{grunbaum}). In this
paper, an algorithm for automatic detection of such a tiling's spiral
structure and its first implementation results will be discussed. Finally, the definitions for several types of spiral tilings will be refined based on this investigation.

\begin{center}
\includegraphics[width=12.5cm]{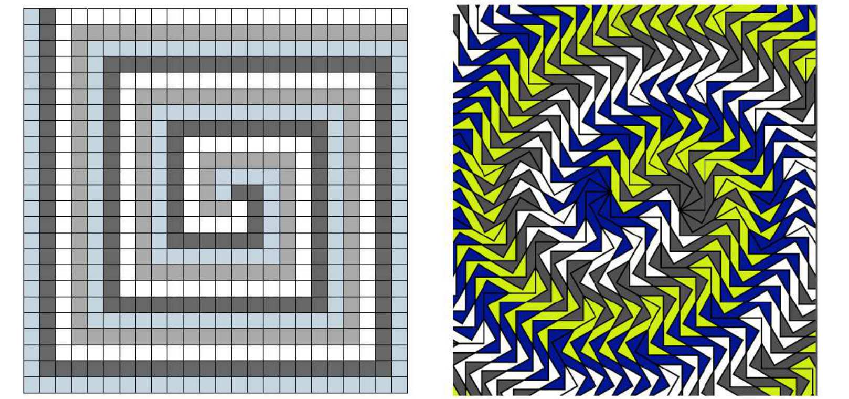} 
\end{center}
\begin{description}
\small
\item [{Figure\ 1}] Spiral structure ``by coloring'' (left) vs. ``by construction'' (right)
\end{description}
If in Figure 1 all colors of tiles were erased and only the tile structure
remained, in the left tiling (of simple squares) nobody could
find a spiral character. On the other hand, the right tiling contains
the spiral structure by construction, although not so easily recognized
without coloring.

One key aspect of the definition of spiral tilings in \cite{klaassen} is that it gives a basis to distinguish between tilings in which the spirals were just introduced by coloring
and those which incorporate a spiral structure.
The two tilings of Figure 1 represent both types of spirals. For the latter type we will define a 
further distinction at the end of this paper.

During this study, we assume that all tiles are closed topological disks. If not specified explicitly, we assume that no singular points exist, where the tiles are clustered. All investigated tilings without such singular points are assumed to be $k$-hedral, which means that there are only finitely many congruence classes. 
We will refer to the definitions from \cite{klaassen} throughout this paper,
so, we decided to put them into the appendix to have them at hand.
First we need the term \textit{L-tiling} which can be summarized using
ordinary language:

``An \textit{L-tiling} allows a partitioning into several parts (called
\textit{arms}), in each of which we can draw a continuous, unlimited
curve (called \textit{thread}) running through the interior of each tile (of the part)
exactly once and winding infinitely often around a certain point''.
In the appendix the reader may have a look for this definition in
strict mathematical terms, but for the further understanding this
one-sentence-version should be sufficient. (The left hand part of
Figure 1 serves as an example of an L-tiling.)

Also the term \textit{S-tiling} from \cite{klaassen} can be summarized in
ordinary words: ``An\textit{ S-tiling} must have the properties of
an L-tiling plus an extra property that neighboring tiles within each
arm are positioned to each other in a way that cannot occur with
two neighbored tiles from different arms (except at the beginning
of an arm)''. E.g., a closer inspection of Figure 1 (right) shows
that within each arm (equal color) there are just two different constellations
of neighboring tiles, and both constellations do not occur with tiles
of different colors. So, this is an S-tiling. (See again the appendix
for a more rigid definition.)

Then for our algorithm we need another pair of definitions.

\begin{spacing}{0}
\noindent \begin{center}
\includegraphics[width=10cm]{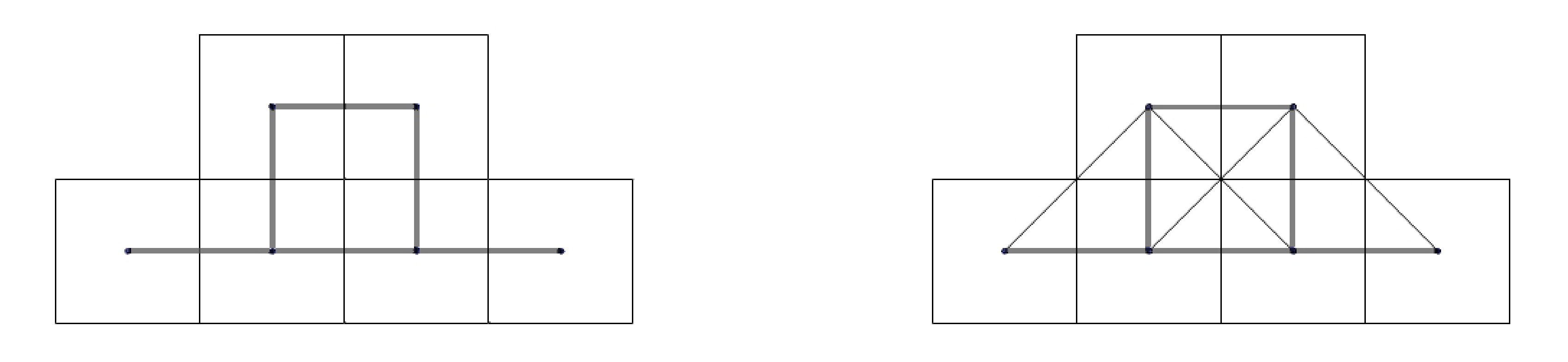}
\par\end{center}
\end{spacing}
\begin{description}
\small
\item [{Figure\ 2}] 6-square subset \textit{M} with \textit{DG(M)} in
grey (left) and \textit{CG(M)} (right)
\end{description}

\begin{description}
\item [{Definition}] \textbf{\textit{direct contact graph (DG)}}\textbf{ and }\textbf{\textit{
contact graph (CG)}}\textbf{:} 
\end{description}
Let $M$ be a connected set of some tiles of a tiling. Then the \textit{``contact
graph of}$M$'' or \textbf{\textit{CG(M)}} is the graph in which each
tile of $M$ is represented by a node and two of such nodes are connected
by an edge iff the corresponding tiles ``have contact'', i.e., have
non-empty intersection \cite{buchsbaum}. We can construct a subgraph of \textit{CG(M)}
called ``\textit{direct contact graph of} $M$'' or \textbf{\textit{DG(M)}}
by deleting all edges for each pair of tiles which share only a finite
number of points of their boundaries. (In graph theory this would be called \textit{dual graph} where the 
tiling is interpreted as a planar graph.)

Figure 2 shows a simple example of \textit{DG} and \textit{CG} for
a small subset of the square tiling. Observe that \textit{CG} in many
cases will not be planar. Both\textit{ DG} and \textit{CG} can be
finite or infinite, depending on the choice of \textit{M}.

\section{The algorithm}

Looking at the above-mentioned definitions for S-tilings, we
observe that they start with a partitioning of the tiling, but do
not tell us how to find it (in our example in Figure 1 ``partition''
and ``coloring'' are equivalent). So, if any automatic recognition
is possible, it must deliver a partitioning into ``spiral
arms''. Given these partitions (or arms in terms of our definition)
it is clear how to proceed further: Check whether a continuous curve
(a so-called thread) can be found satisfying the necessary conditions.
For practical reasons, we decided to search for Hamilton paths \cite{harel} 
within each candidate for a spiral arm. Although this is not  equivalent
to definition L or S, for a huge subset of S-tilings (maybe for all of them) 
the spiral arms can be regarded as Hamiltonian w.r.t \textit{DG} or\textit{ CG}.  This  can
be easily implemented using graph libraries. (A Hamilton path within a connected component of
\textit{DG} or\textit{ CG} means that we can walk through the component
along its edges meeting every node exactly once.) 

Let us first describe the main ideas of the algorithm just by words:
\begin{itemize}
\item Build classes of neighboring tile pairs according to their relative
position to each other
\item For each possible subset of these classes cut the tiling at the intersection
of each tile pair belonging to one of the selected classes
\item After each cut check the resulting connected components whether they
allow a Hamilton path running through each component of the (direct) contact graph winding around a central point 
\end{itemize}

To give an example, in Figure 2.1 we can find four classes of tile
pairs (one connected by a short edge and three others sharing a long
edge in different ways). 

\begin{frame}{}
\hbox{\hspace{-0.5em} \includegraphics[width=8cm]{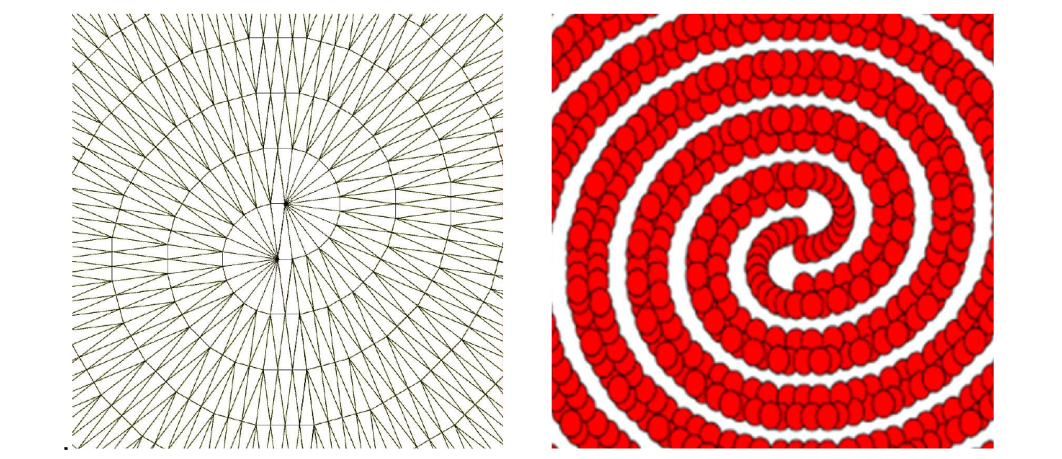} }
\end{frame}
\begin{description}
\small
\item [{Figure\ 2.1}] Example: spiral tiling (left) and result from the algorithm (right)
\end{description}

It is obvious that just the connections via
``short'' edges had to be deleted by the algorithm to find the spiral
structure. In this case it is a two-arm spiral where both arms meet
at the center. We see (at the right half of Figure 2.1) the direct
contact graph (\textit{DG}) of the tiling after the described cut. 

Now the more formal description must follow: For an algorithmic approach,
we have to restrict our scope to finite portions of a given tiling.
Throughout this section let $M$ be the investigated portion of a
tiling that should be checked by our algorithm. In section 4 the appropriate
choice of\textit{ M} will be addressed. Those tiles in $M$ which
are (in the unlimited tiling) neighbors of tiles lying outside of
$M$ are called the \textit{border} $B_{M}$. Then \textit{CG(M)}
(or \textit{CG} for short) is the contact graph of \textit{M} and\textit{
DG(M)} (or \textit{DG}) the corresponding direct contact graph. We
start with classifying all edges in these graphs depending on how
their corresponding tile pairs in $M$ are positioned to each other: 

For each edge $k$ of $DG$ we form the class {[}$k${]} consisting
of all the edges $k'$ of $DG$ for which the two tiles (determined
by the endpoints of $k'$) are congruent to the corresponding
tiles determined by $k$, through an orientation-preserving isometry
of the plane (that is, by translation or rotation). Let the set of
all such classes be denoted $K_{M}=\left\{ [k_{1}],[k_{2}],\text{\dots}\right\} $.\footnote{One could call such classes \emph{edge classes} or \emph{tile pair classes}, which is equivalent here.}
(During the algorithm we will also need additional edge classes from
$CG$, constructed in the same way.) For a class {[}$k${]} we consider
the set $E(k)$ of edges in {[}$k${]}. For
each subset of classes \textbf{\textit{K}} $\subset K_{M}$ we write
\textit{E(}\textbf{\textit{K}}\textit{)} 
for the corresponding edge set as a union of $E(k)$. 
An edge from \textit{E(}\textbf{\textit{K}}\textit{)}  between the tiles $T_1$ and $T_2$ should be denoted as $(T_1,T_2)$.
For each of these subsets of edge classes
in $DG$ we can check what happens if all these edges were deleted.
How do the remaining connected components of $DG$ ``behave''? Several
steps were included in order to exit the loops as early as possible.

For shortness, we will use the term ``component''
for ``connected component''. 
\begin{description}
\item [{Algorithm:}] First check whether there are at least three congruent
tiles differing in orientation or reflection within $M$.\\
If not, end the algorithm with empty result.\\
Else, form the set of classes $K_{M}$ as described above. \\
Next we define an operation to be performed on each nonempty \textbf{\textit{K}}$\subset$\,$K_{M}$.
\item [{Operation\ A:}] (using \textbf{\textit{K }}as input and returning
either \textbf{\textit{K }}plus a graph or the result ``discarded''
if \textbf{\textit{K}} cannot fulfill a condition)\end{description}
\begin{verse}
Check whether all components of $\bigcup\limits_ {(T_i,T_j)\in E(\boldsymbol{K})} \!\!\!\!\!\!\!\!\! T_i \cap T_j \,\,\,$
are connected to the border $B_{M}$; [Remark: Represents boundaries of spiral arms.] 

if not, discard \textbf{\textit{K}} and finish Operation A with result
``discarded''; 

if yes, delete the edges \textit{E(}\textbf{\textit{K}}\textit{)}
from $DG$, the result is called $G$; 

if all components of $G$ are connected to $B_{M}$ and allow a Hamilton path without self-intersections, go directly
to ({*});

else, do the following steps with \textbf{\textit{K}} plus any combination
of edge classes from $CG$, called ``\textbf{\textit{K}}-extension'', each of which generates a new $G$: 

If for such an extended \textbf{\textit{K}} a component of the new
$G$ is not connected to $B_{M}$ or does not allow a self-avoiding Hamilton path, ignore this \textbf{\textit{K}}-extension and try the next possible one;

if a tile in $M$ has more than two vertices where it meets other tiles at single points (being connected to these tiles by edges in the new $G$), ignore this \textbf{\textit{K}}-extension and try the next possible one; [Remark: Excluding cases like the checkers tiling in [1] Figure 3, where a spiral arm is not simply connected.] 
 
({*})  

If for each component of $G$ the number of tile equivalence classes w.r.t translation is less than 3 \textendash{} discard (extended) \textbf{\textit{K}} and continue with the next one (if
extensions were needed); 

if extensions were needed, return all non-discarded variants of \textbf{\textit{K}}
plus $G$ or ``discarded'' as result when all extensions were checked; 

else return (\textbf{\textit{K}}, $G$) if non-discarded or else
return ``discarded''. 

(\textbf{End of A})

Perform Operation A with all nonempty subsets of $K_{M}$. All non-discarded
subsets are candidates for spiral partitions. Sort the non-discarded
subsets by the number of components of $G$ (= nbr. of arms) in increasing
order.\end{verse}
\begin{description}
\item [{Operation\ B:}] (using each non-discarded $G$ as input if there
is any)\end{description}
\begin{verse}
For each component of $G$: Find a continuous piece-wise linear
curve through the corresponding tiles following the possible Hamilton
paths and modify it to check whether the conditions for being a thread
can be fulfilled. {[}Remark: this section of the algorithm is not
difficult for the human eye but needs considerable programming efforts.
On the other hand, by methods of computer graphics and optimization
this task could be handled in principle. Since the above-mentioned
``psychological effect'' is not needed here, we decided not to code
this section of the algorithm.{]} 

If one component does not allow a thread, $G$ has to be discarded. As a final result the components of $G$ each with a valid thread
represent the spiral arms.
\end{verse}
\textbf{End of Algorithm}

\ 

This algorithm (except for Operation B) was implemented in Python,
which is by far not the fastest language but offers a lot
of libraries for graph operations. 

Let us return to the example in Figure 2.1. The separation into two
arms cannot be managed by the implemented algorithm, but the spiral structure
was recognized. Only the direct contact graph is needed in this case,
but there will be some examples with \textit{CG} in the following
section.

\section{Results\ }

As a set of test cases we took several tilings from \cite{wichmann} with spiral structure. 

\begin{spacing}{0}
\noindent \begin{center}
\includegraphics[width=12.2cm]{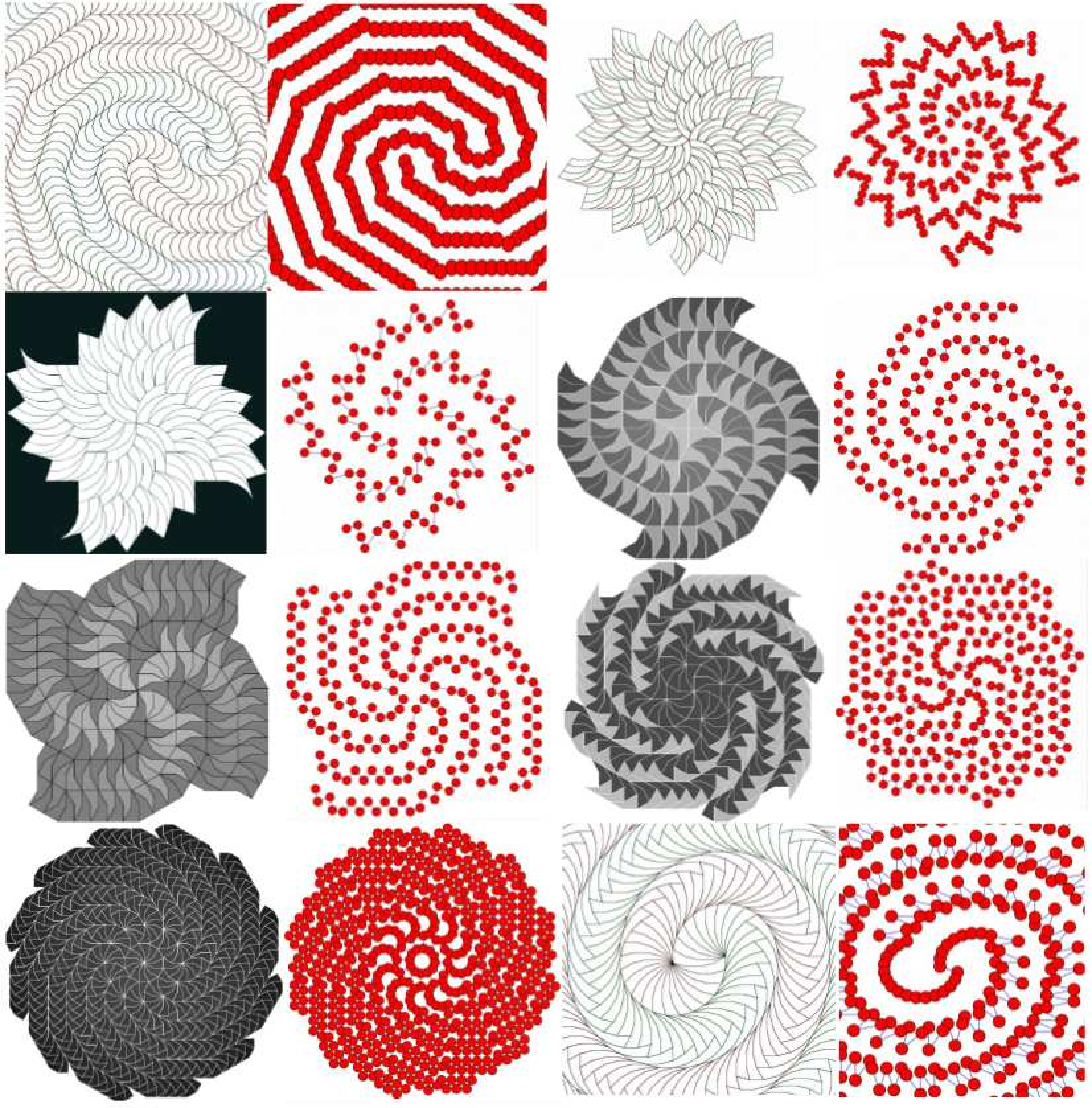}
\par\end{center}
\end{spacing}
\begin{description}
\small
\item [{Figure\ 3.1a}] Tilings and resulting graphs from the algorithm 
\end{description}
In the first and third column of Figure
3.1a and 3.1b we show the tilings and right hand besides them in the
second and fourth column the resulting graphs from our algorithm.
For the majority of tilings the algorithm works with \textit{DG}.
The list of examples is continued with Figure 3.1b (still with
usage of \textit{DG} instead of \textit{CG}).

\begin{spacing}{0}
\noindent \begin{center}
\includegraphics[width=12.2cm]{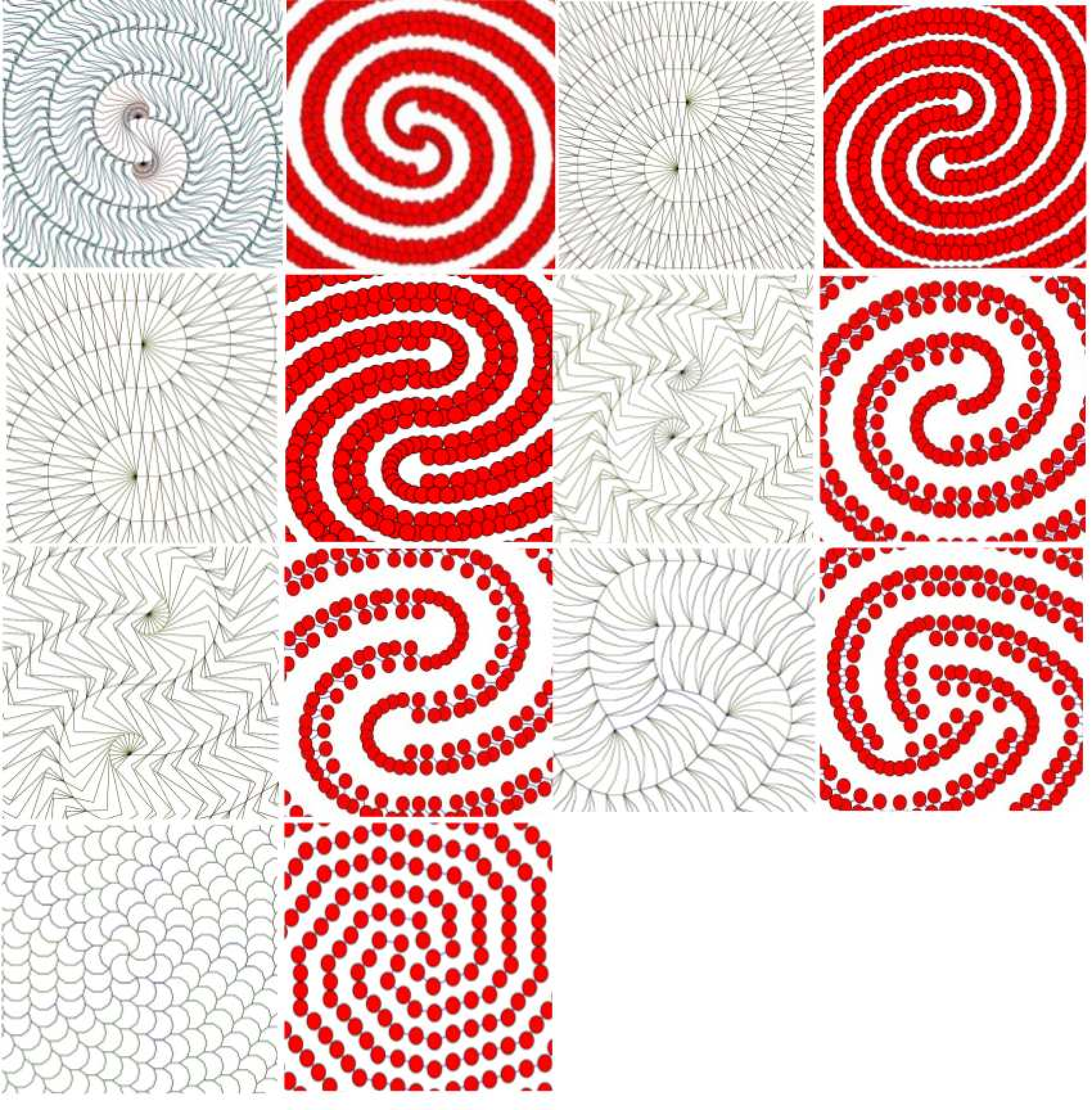}
\par\end{center}
\end{spacing}
\begin{description}
\small
\item [{Figure\ 3.1b}] Tilings and resulting graphs from the algorithm 
\end{description}

\begin{spacing}{0}
\noindent \begin{flushleft}
\includegraphics[width=12.2cm]{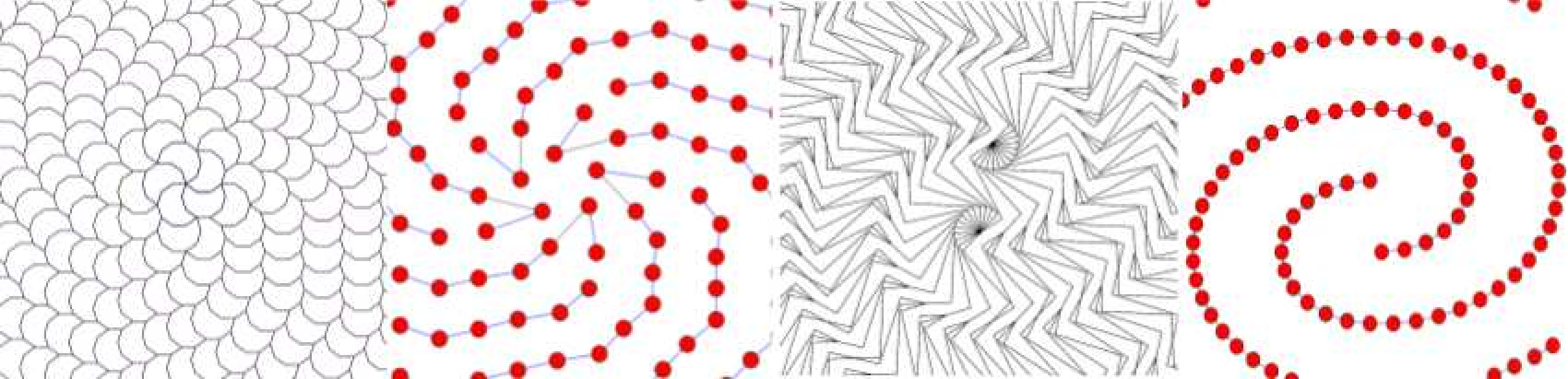}
\par\end{flushleft}
\end{spacing}
\begin{description}
\small
\item [{Figure\ 3.2}] Tilings and resulting graphs with usage of \textit{CG} 
\end{description}
There are some rare cases where \textit{CG} is needed (see Figure 3.2). So, it is recommended to start with \textit{DG} and only
if nothing could be found, a second round with \textit{CG} should
be performed.

It is interesting to note that the algorithm also works for one-armed
spirals. Though the definition for this type differs slightly from
Definition S (see Appendix: Definition O for more details) the algorithm
(up to Operation B) can be applied without any changes. Operation B can be performed here in simplified version, since only the spiral boundary curve has to be checked if it is winding around its starting point. In the above-mentioned collection
of spiral tilings \cite{wichmann} there are two examples for this case (Figure
3.3).

\begin{spacing}{0}
\noindent \begin{flushleft}
\includegraphics[width=12.2cm]{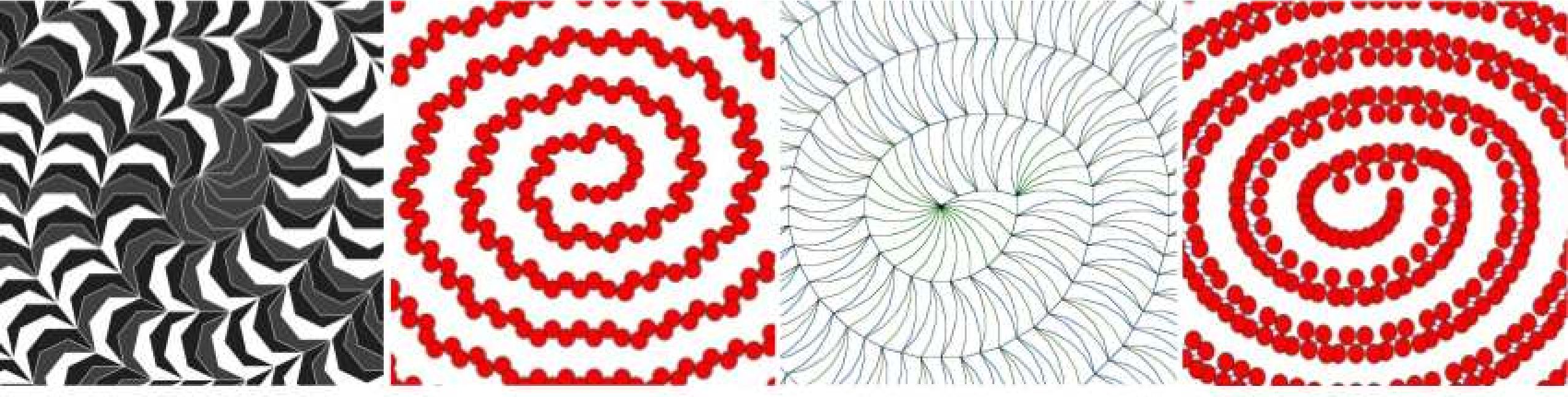}
\par\end{flushleft}
\end{spacing}
\begin{description}
\small
\item [{Figure\ 3.3}] Tilings and resulting graphs for the one-armed case 
\end{description}
There are some special situations, where the results indicate more
than one spiral partitioning. In Figure 3.4 we show two different
spirals for the same tiling that were both found by the algorithm

\includegraphics[width=11cm]{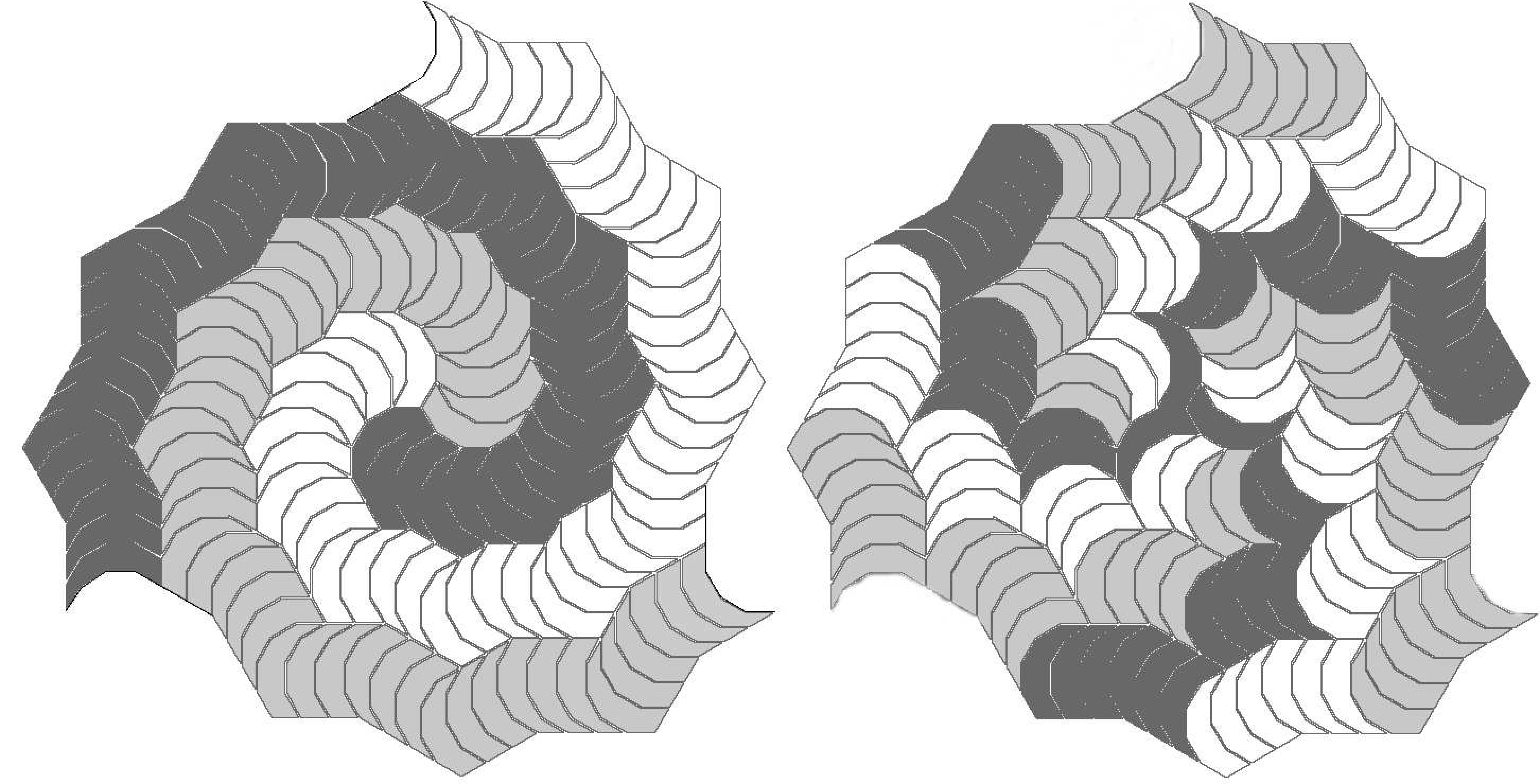}
\begin{description}
\small
\item [{Figure\ 3.4}] Two different S-partitions for the same tiling 
\end{description}
The spiral arms on the right side of Figure 3.4 do not look very ``natural''
compared to the spirals on the left half, but they fulfill all conditions
for an S-tiling. Only the heuristic argument could be applied that
the partition with lower number of arms should be preferred. This
is the reason for the final part of the algorithm, where a sorting
of the resulting graphs has been proposed to find the result with
lowest number of connected components.

\section{Complexity and other algorithmic aspects}

It is quite obvious that an algorithm containing a loop over all subsets
of a given finite set must have exponential complexity (w.r.t. the
number of edge classes $\vert K_{M}\vert$). Hence, there will be
cases where the algorithm's runtime outruns all practical limits.
It should be noted here that the whole investigation did not aim on
efficient implementation, but to answer the question whether such
an algorithm exists at all. 

What can be done now in cases of extremely long runtime? Such examples
exist, but we are lucky that they are rare. For these few cases we
propose to apply an algorithmic test ``by hand'' in a way that the
following items should be checked to decide whether the algorithm
will (or won't) be successful. We use again the simple structure of
Figure 2.1 to illustrate the steps:
\begin{itemize}
\item classify all edges of the direct contact graph\textit{ DG} by assigning
integers for each class of direct neighbors to define the edge classes
$K_{M}$ (In Figure 2.1 there are four classes: Let us assign 1 to
the neighbors sharing a short edge and 2, 3 and 4 to the other classes
of neighbors sharing a long edge.) 
\item if spiral arms can be observed by the human eye:
Consider the spiral arms as subsets of $M$ and
 run along their boundaries
 to find the specific subset \textbf{\textit{K}} of $K_{M}$.
If a spiral arm locally shrinks to a single point, as in
Figure 1 (right), go back to the previous item but use \textit{CG}
(In our example in Figure 2.1 just \textit{DG} is needed and the arms'
boundaries are easily characterized just by the short edges, so we
choose \textbf{\textit{K}} = \{1\}.) 
\item check whether the chosen \textbf{\textit{K}} finishes Operation A
without being discarded (This is easily checked in our example since
the tiling contains more than three tiles with different rotation
angles and each component of the resulting $G$ - after deleting
the connections via short edges - can be naturally traversed by a
Hamilton path. All these paths are connected to the outside border
region, which is also true for the arms' boundary.)
\item perform Operation B for the components of the non-discarded results
of Operation A, i.e., find a thread - or maybe several of them - following
the Hamilton path(s) (In our example this is done straight forward
with two threads starting close to the tiling's center.)
\end{itemize}
If by these checks a single subset\textbf{\textit{ K}} is found not
being discarded by Operation A, it is shown that the algorithm must
find this result within finite time. All remaining tilings from the
literature (less than 10) were investigated with the result that in
all cases where definition S is satisfied, the algorithm will return
a spiral partition. Also the somehow unexpected spiral structure within
the Hirschhorn tiling  can be detected by this analysis (discussed in the last section, see Figure 5.2).
 In the same section we will see another simple example which demonstrates the
advantages of the algorithm's application ``by hand''. 

There is one interesting case in Brian Wichman's collection \cite{wichmann}
showing kind of disrupted spiral arms (see Figure 4.1). The two arms
indicated by two different colors are following a spiral structure
from inside to outside, but it is not possible to draw a continuous
path following the spiral within the interior of each arm. 

\ 

\begin{center}
\includegraphics[scale=0.4, angle=90]{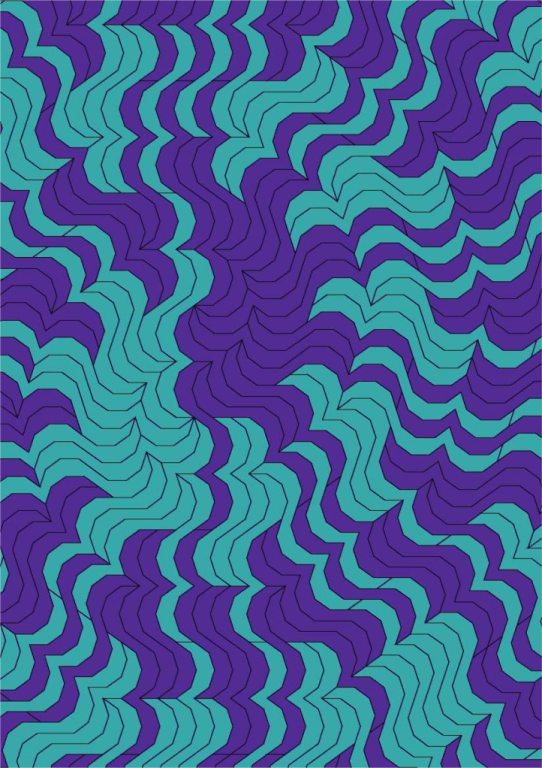}
\end{center}
\begin{description}
\small
\item [{Figure\ 4.1}] A tiling from \cite{wichmann} with disrupted spiral arms 
\end{description}
Our algorithm (here not by hand but by software) returned a negative
result in this case, which is correct since neither definition S nor
even L can be satisfied. 

By its nature, an algorithm working on a finite portion of the tiling
cannot in all cases distinguish between ``true'' spirals and partitions
which start like spirals but later stop the spiral behavior. (Figure
11 in \cite{klaassen} shows an example of such a pseudo-spiral partition.)
Therefore one could start the algorithm with a smaller part of the tiling as described 
and then add further tiles outside of the so-called border region.
Then one could check whether the orientation of the tiles within a
spiral arm will further change or remains in one or two fixed angular
positions. In addition, a further difficulty could occur: It is not sure that
the spiral center is always placed in the middle of the finite portion
of the tiling. So, one might first look for this center by searching
for the part where the highest number of tiles with different orientation
are clustered.

The proposed refinements from this section are all possible in principle
but the described version of the algorithm worked well enough without
it.

\section{Discussion and further refinements}

The main result of this paper is the fact that an algorithm can be
designed to decide whether a given tiling has or doesn't have a spiral
structure. This is done by a method of partitioning into spiral arms.
As we have seen in the results section, the proposed algorithm can
be applied to a wide range of tilings. We can claim that all known
spiral tilings from the literature (in the meaning of definition S
or O resp.) can be detected by the algorithm. The vast majority was
covered by our Python implementation while the remaining part (less
than 10) could be analyzed ``by hand'' following the algorithmic
check list described in the previous section. So, the algorithm is
working as desired with the limitation of not being very efficient
for all cases due to its exponential complexity.

This application ``by hand'' can also be used to decide
whether a given tiling contains more than one spiral structure.

\begin{spacing}{0}
\noindent \begin{center}
\includegraphics[width=10cm]{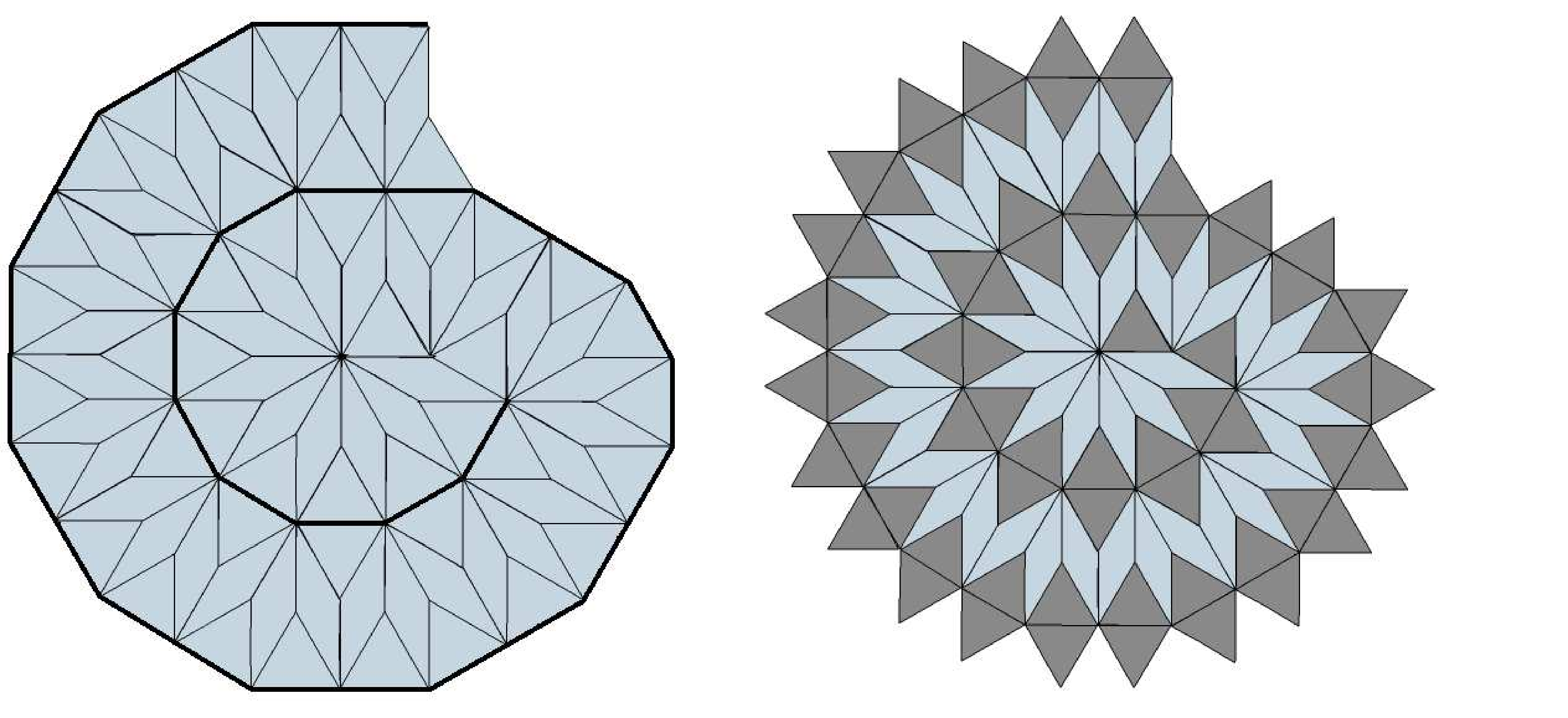}
\par\end{center}
\end{spacing}
\begin{description}
\small
\item [{Figure\ 5.1}] The algorithm applied ``by hand'' to analyze a tiling 
\end{description}

We can demonstrate this with a tiling presented in \cite{klaassen} to find
out whether more than one spiral arm exists in this case (see Figure
5.1). First, we cut the tiling at those edges shared by two tringles
(= thick line) to get the left hand version (one spiral arm). Alternatively
- on the right hand side - we cut the same tiling at those edges where
a triangle meets a rhombus (= edges shared between dark grey and light
grey tiles) to find the right hand partitioning (two spiral arms).
This means that the algorithm ``by hand'' can also be used as a
method to partition a given tiling for a better understanding of its
structure. 

The concept of structure analysis developed in this paper can be used for other structures than spirals, as well.  In any case the final result is a tile set partition, where in each part the tiles are positioned to each other in a different way than on the parts' boundaries. So, one could ask, what typical structures could be found in this way: For the large domain of periodic tilings, we will often find partitions in form of stripes or patches. For non-periodic tilings, especially with rotational symmetry  - but not restricted to those -  we will detect ring-like structures, where each ring is surrounded by a larger one. For such ring partitions, we distinguish two types which can be defined in a way analogous to definition L and S.

\begin{description}
\item [{Definition}] \textbf{\textit{weak ring partition}: }
A tile set partition of a plane tiling into infinitely many parts is called a \textit{weak ring partition} if each part (as union of its tiles) 
contains a closed Jordan curve $\theta$ (called \textit{thread}) around a fixed central point,
$\theta(t)=r(t)\exp{(i\varphi(t))}$ with the plane identified with $\ensuremath{\mathbb{C}}$, 
 $r(t)>0$, $t\in [0,1]$ and $\varphi$ being monotonic with $\varphi([0,1])= [0,2\pi]$.
 For each tile $T$ in the part the intersection of the interior
of $T$ with the image of $\theta$ is nonempty and connected. 
The threads do not meet or cross each other.
\end{description}

Note that this definition could also be used for tilings with a singular point, where arbitrarily small tiles are clusterd. Apart from this, a huge number of tilings allow weak ring partitions, however, it is not a simple question how to characterize the family of tilings that share this property. We can pose this as an open problem so far.

For the further analysis we need a stronger version of this definition. The condition is analogous to S2 from definition S with `arm' replaced by `part':

\begin{description}
\item [{Definition}] \textbf{\textit{strong ring partition}: }
A tile set partition of a plane tiling with all properties of a weak ring partition is called a \textit{strong ring partition} if an additional condition holds:
If any two tiles $T_{1},T_{2}$ in a part are
direct neighbors and can be respectively
mapped by an operation $\tau$ (composed of translation, rotation or scaling) onto another tile pair $\tau(T_{1})$
and $\tau(T_{2})$, these must also be direct neighbors within
a part.  ($T_{1},T_{2}$ from the same part are called
\textit{direct neighbors} if $T_{1}\text{\ensuremath{\cap}}T_{2}$
is cut\footnote{Here and in all other occurrences "cut by the thread" means that the thread (by passing from $T_{1}$ to $T_{2}$)  intersects  $T_{1}\text{\ensuremath{\cap}}T_{2}$, which might also be just a single vertex.} by the part's thread or contains more than a finite number
of points.) 
\end{description}

Figure 5.2 (left) shows an example of a strong ring partition. The scaling operation was inserted here to make this definition applicable also in the context of tilings with singular points, see below in this section. It is obvious that by the techniques of the algorithm presented in section 2 one could automatically check whether a tiling allows or doesn't allow a strong ring partition. 

Now we can separate tilings with a spiral structure from those with a ring structure, which sometimes can both occur simultaneously.

\begin{spacing}{0}
\noindent \begin{center}
\includegraphics[width=11cm]{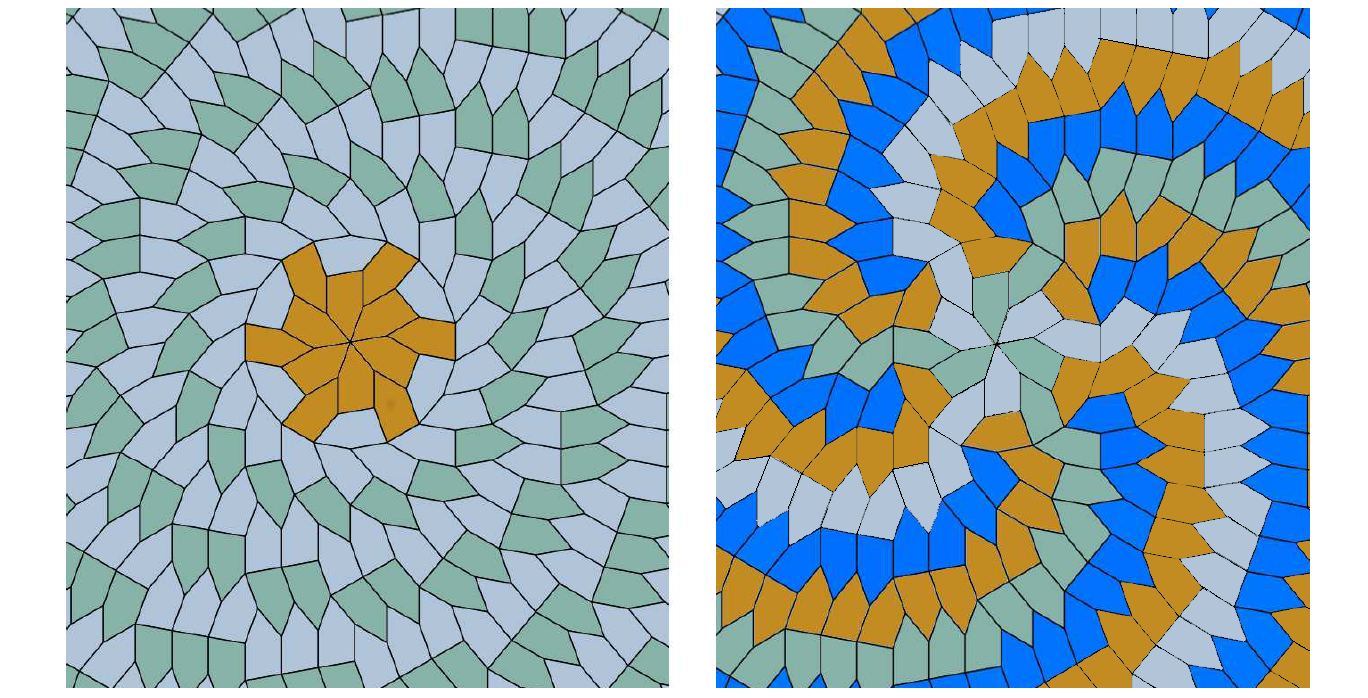}
\par\end{center}
\end{spacing}
\begin{description}
\small
\item [{Figure\ 5.2}] A tiling with strong ring partition (left) and S-partition (right)
\end{description}

\begin{description}
\item [{Definition:}] 
A $k$-hedral tiling is called a \textit{strong spiral tiling} (respectively \textit{strong S-tiling} or \textit{strong O-tiling})
if it is an S- or O-tiling and additionally \textit{doesn't} allow a strong ring partition. (Hence, strong spiral tilings and strong ring partitions exclude each other.) 
\end{description}

A closer inspection shows that most of the known S- or O-tilings are also strong S- respectively strong O-tilings. A famous example where this is \textit{not} the case is the Hirschhorn tiling. In Figure 5.2 we can observe on the left side the obvious ring structure and on the right side the spiral arms allowing an S-partition. 

In the context of tilings with one singular point, we can do the same to separate the ring structure from the spirals.   

\begin{description}
\item [{Definition:}] 
A tiling with one singular point and finitely many similarity classes is called a \textit{strong spiral tiling} (or \textit{strong P-tiling})
if it allows a partition according to definition P but \textit{doesn't} allow a strong ring partition.  
\end{description}

Tilings with strong ring partition and spiral structure - regardless of having a singular point or not - as shown in Figure 5.2 or in \cite{staana} often have the property that the spirals are in some sense hidden or visually dominated by the ring structure. They can be viewed as ``picture puzzles''. So, though we have demonstrated that spirals (and other structures) can be detected principally without human aid by algorithms, in the context of perception the quote from the beginning remains true that ``to some
extent, at least, the spiral effect is psychological\textquotedblright.

\section*{Acknowledgements}

Although it is no longer possible in person, I would like to express
my gratitude to the late Branko Gr\"unbaum for the fruitful discussion
during the development of this paper and thanks to Brian A. Wichmann
for his support with the tilings from his great collection.

\section*{Appendix: Definitions from \cite{klaassen} }
\small
Here are the definitions from \cite{klaassen} which are used in this paper.
In section 2 we tried to paraphrase them using ordinary language. 
\begin{description}
\item [{Definition\ L}] (spiral-like) :
\end{description}
A partition of a plane tiling into more than one separate classes
(called \textit{arms} here) is defined as a \textit{spiral-like partition}
or \textit{L-partition} under the following conditions. (The plane
is identified with the complex plane $\ensuremath{\mathbb{C}}$ where
the origin is represented by a selected point of the tiling.) 
\begin{description}
\item [{L1:}] \, For each arm $A$ (as a union of tiles from one class) there
exists a curve $\theta:\mathbb{R}_{0}^{+}\rightarrow A\text{\,}$$\subset\mathbb{C}$
with $\theta(t)=r(t)\exp{(i\varphi(t))}$ called a \textit{thread}
, where both $r$ and $\varphi$ are continuous and unbounded and
$\varphi$ is monotone. Curve $\theta$ does not meet or cross itself
or any thread from another arm of the tiling. 
\item [{L2:}] \, For each tile $T$ in $A$ the intersection of the interior
of $T$ with the image of $\theta$ is nonempty and connected. 
\item [{Remark:}] A plane tiling with an L-partition shall be called an
\textit{L-tiling} or a \textit{spiral-like tiling}. 
\end{description}
\ 
\begin{description}
\item [{Definition\ S:}] (for tilings with more than one spiral arm)
\end{description}
A partition of a plane tiling is defined as a \textit{spiral partition}
or \textit{S-partition} under the following conditions. 
\begin{description}
\item [{S1:}] \, It must be an L-partition (see definition L). 
\item [{S2:}] \, If any two tiles $T_{1},T_{2}\text{\ensuremath{\in}}A$ are
direct neighbors and can be respectively
mapped by a direct isometry $\tau$ onto another pair of tiles $\tau(T_{1})$
and $\tau(T_{2})$, these must also be direct neighbors within
an arm. This rule can be ignored if the image pair contains the beginning
of an arm, i.e., contains $\theta(0)$. ($T_{1},T_{2}\in A$ are called
\textit{direct neighbors} if $T_{1}\text{\ensuremath{\cap}}T_{2}$
is cut by the thread of $A$ or contains more than a finite number
of points.) 
\item [{Remark:}] A plane tiling which allows an S-partition shall be called
an \textit{S-tiling}.
\end{description}
\ 
\begin{description}
\item [{Definition\ O}] (for one-armed spirals)\ :
\end{description}
A tiling of the plane is called a \textit{spiral tiling with one arm} or an \textit{O-tiling}
under the following conditions (The plane is identified with the complex
plane $\ensuremath{\mathbb{C}}$ where the origin is represented by
a selected point of the tiling.)
\begin{description}
\item [{O1:}] \, There exists a curve $b:\mathbb{R}_{0}^{+}\rightarrow\mathbb{C}$
with $b(t)=r(t)\exp{(i\varphi(t))}$ called \textit{spiral boundary},
where both $r$ and $\varphi$ are continuous and unbounded. Curve
$b$ does not meet or cross itself and runs completely on boundaries
of tiles 
\item [{O2:}] \, If $T_{1},T_{2}$ are direct neighbors and can be respectively
mapped by a direct isometry $\tau$ onto another pair of tiles $\tau(T_{1})$
and $\tau(T_{2})$, these tiles
must also be direct neighbors. This rule can be ignored if the image
pair lies at the beginning of the boundary (i.e., contains $b(0)$).
(`\textit{Direct neighbors}' means here that $T_{1}\cap T_{2}$
contains more than a finite number of points but not from the spiral
boundary.) \end{description}

\begin{list}{}%
{\leftmargin=2.5em
\itemindent=0em}
\item [\bf{Definition\,\,P}] (for tilings with one singular point)\,: A tiling
of the Euclidean plane with exactly one singular point together with
a partition is called a \emph{spiral P-tiling}
under the following conditions. (The plane is identified with $\mathbb{C}$
where the origin is represented by the singular point.)
\item [\bf{P1:}]  The partition fulfills L1 and L2 but with $\theta:\mathbb{R}\rightarrow A\text{\,}\subset\mathbb{C}$
and with $\varphi$ being unbounded in both directions.
\item [\bf{P2:}] If any two tiles $T_{1},T_{2}\text{\ensuremath{\in}}A$
are direct neighbors and can be respectively
mapped by an operation $\tau$ (composed of translation, rotation or scaling) onto another pair of tiles $\tau(T_{1})$
and $\tau(T_{2})$, these tiles must also be direct neighbors
within an arm. (`\emph{Direct neighbors}' means here that $T_{1}\text{\ensuremath{\cap}}T_{2}$
is cut by the thread of $A$.)
\end{list}
\end{document}